\documentclass[10pt, a4paper]{article}

\usepackage{amsmath}
\usepackage{amssymb}
\usepackage{amsthm}
\usepackage{mathrsfs}
\usepackage{hhline}
\usepackage{array}
\usepackage{multirow}
\usepackage{cite}
\usepackage{enumerate}
\usepackage{xcolor}
\usepackage{graphics}
\usepackage{epsfig}
\usepackage{subfigure}
\usepackage{longtable}

\usepackage[mathlines]{lineno}
\modulolinenumbers[2]
\newcommand*\patchAmsMathEnvironmentForLineno[1]{%
\expandafter\let\csname old#1\expandafter\endcsname\csname #1\endcsname  \expandafter\let\csname oldend#1\expandafter\endcsname\csname end#1\endcsname  \renewenvironment{#1}%
{\linenomath\csname old#1\endcsname}%
{\csname oldend#1\endcsname\endlinenomath}}%
\newcommand*\patchBothAmsMathEnvironmentsForLineno[1]{%
\patchAmsMathEnvironmentForLineno{#1}%
\patchAmsMathEnvironmentForLineno{#1*}}%
\AtBeginDocument{%
\patchBothAmsMathEnvironmentsForLineno{equation}%
\patchBothAmsMathEnvironmentsForLineno{align}%
\patchBothAmsMathEnvironmentsForLineno{flalign}%
\patchBothAmsMathEnvironmentsForLineno{alignat}%
\patchBothAmsMathEnvironmentsForLineno{gather}%
\patchBothAmsMathEnvironmentsForLineno{multline}%
}

\definecolor{vividviolet}{rgb}{0.62, 0.0, 1.0}
\def\rsq{\hspace*{\fill}$\blacksquare$\medskip}
\def\rs{\mathcal R}
\def\cs{\mathcal C}
\def\d{\mathcal D}
\def\M{\mathcal M}
\def\S{\mathcal S}
\newtheoremstyle{rem}
  {10pt}          
  {10pt}  
  {\rm}  
  {}
  {\bf}  
  {: }    
  { }    
  {}     
\theoremstyle{rem}
\newtheorem{rem}{Remark}[section]
\newtheorem{example}{Example}[section]
\newtheorem{problem}{Problem}[section]

\newtheoremstyle{theorem}
  {10pt}          
  {10pt}  
  {\it}  
  {}
  {\bf}  
  {: }    
  { }    
  {}     
\theoremstyle{theorem}

\newtheorem{theorem}{Theorem}[section]
\newtheorem{lemma}[theorem]{Lemma}

\newtheorem{corollary}[theorem]{Corollary}

\numberwithin{equation}{section}
\def\Z{\mathbb{Z}}
\def\N{\mathbb{N}}
\newtheorem{defi}{Definition}[section]

\usepackage{graphicx}
\graphicspath{%
    {converted_graphics/}
    {/}
}

\def\ms{\medskip}
\def\nt{\noindent}

\begin{document}
\baselineskip18truept
\normalsize
\begin{center}
{\mathversion{bold}\Large \bf On local antimagic total labeling of amalgamation graphs }

\bigskip
{\large  Gee-Choon Lau{$^{a,}$}\footnote{Corresponding author.}, Wai-Chee Shiu{$^{b}$}}\\

\medskip

\emph{{$^a$}Faculty of Computer \& Mathematical Sciences,}\\
\emph{Universiti Teknologi MARA (Segamat Campus),}\\
\emph{85000, Johor, Malaysia.}\\
\emph{geeclau@yahoo.com}\\

\medskip

\emph{{$^b$}Department of Mathematics,\\ The Chinese University of Hong Kong,}\\
\emph{Shatin, Hong Kong.}\\

\emph{wcshiu@associate.hkbu.edu.hk}\\

\end{center}

\begin{abstract}
Let $G = (V,E)$ be a connected simple graph of order $p$ and size $q$. A graph $G$ is called local antimagic (total) if $G$ admits a local antimagic (total) labeling. A bijection $g : E \to \{1,2,\ldots,q\}$ is called a local antimagic labeling of $G$ if for any two adjacent vertices $u$ and $v$, we have $g^+(u) \ne g^+(v)$, where $g^+(u) = \sum_{e\in E(u)} g(e)$, and $E(u)$ is the set of edges incident to $u$. Similarly, a bijection $f:V(G)\cup E(G)\to \{1,2,\ldots,p+q\}$ is called a local antimagic total labeling of $G$ if for any two adjacent vertices $u$ and $v$, we have $w_f(u)\ne w_f(v)$, where  $w_f(u) = f(u) + \sum_{e\in E(u)} f(e)$. Thus, any local antimagic (total) labeling induces a proper vertex coloring of $G$ if vertex $v$ is assigned the color $g^+(v)$ (respectively, $w_f(u)$). The local antimagic (total) chromatic number, denoted $\chi_{la}(G)$ (respectively $\chi_{lat}(G)$), is the minimum number of induced colors taken over local antimagic (total) labeling of $G$. In this paper, we determined $\chi_{lat}(G)$ where $G$ is the amalgamation of complete graphs.\\ 


\end{abstract}

\section{Introduction}
Consider a $(p,q)$-graph $G=(V,E)$ of order $p$ and size $q$. In this paper, all graphs are simple. For positive integers $a < b$, let $[a,b]=\{a,a+1,\ldots,b\}$. Let $g:E(G)\to [1,q]$ be a bijective edge labeling that induces a vertex labeling $g^+: V(G) \to \N$ such that $g^+(v) = \sum_{uv\in E(G)} g(uv)$. We say $g$ is a {\it local antimagic labeling} of $G$ if $g^+(u)\neq g^+(v)$ for each $uv\in E(G)$~\cite{Arumugam,Bensmail}. The number of distinct colors induced by $g$ is called the {\it color number} of $g$ and is denoted by $c(g)$. The number
\[\chi_{la}(G)=\min\{c(g)\;|\; g \mbox{ is a local antimagic labeling of } G\}\]
is called the {\it local antimagic chromatic number} of $G$~\cite{Arumugam}. Clearly, $\chi_{la}(G)\ge \chi(G)$.

\nt Let $f: V(G)\cup E(G) \to [1,p+q]$ be a bijective total labeling that induces a vertex labeling $w_f : V(G) \to \N$, where
$$w_f(u)=f(u) + \sum_{uv\in E(G)} f(uv)$$
and is called the {\it weight} of $u$ for each vertex $u \in V(G)$. We say $f$ is a {\it local antimagic total labeling} of $G$ (and $G$ is {\it local antimagic total}) if $w_f(u)\ne w_f(v)$ for each $uv\in E(G)$. Clearly, $w_f$ corresponds to a proper vertex coloring of $G$ if each vertex $v$ is assigned the color $w_f(v)$. If no ambiguity, we shall drop the subscript $f$. Let $w(f)$ be the number of distinct vertex weights induced by $f$. The number
\[\min\{w(f)\;|\;f \mbox{ is a local antimagic total labeling of } G\}\]
is called the {\it local antimagic total chromatic number} of $G$, denoted $\chi_{lat}(G)$. Clearly, $\chi_{lat}(G)\ge \chi(G)$. It is well known that determining the chromatic number of a graph $G$ is NP-hard~\cite{Zukerman}. Thus, in general, it is also very difficult to determine $\chi_{la}(G)$ and $\chi_{lat}(G)$.

\nt For a graph $G$, the graph $H = G\vee K_1$ is obtained from $G$ by joining a new vertex to every vertex of $G$.  We refer to~\cite{Bondy} for notation not defined in this paper.

\ms\nt In~\cite{Haslegrave}, the author proved that every connected graph of order at least 3 is local antimagic. In~\cite{Lau+K+S}, the authors proved that every graph is local antimagic total. We shall need the following theorem in~\cite{Lau+K+S}.  

\begin{theorem}\label{thm-K1VG} Let $G$ be a graph of order $p\ge 2$ and size $q$ with $V(G)=\{v_i\,|\,1\le i\le p\}$. 

\begin{enumerate}[(a)]
  \item $\chi(G)\le \chi_{lat}(G) \le \chi_{la}(G \vee K_1) - 1$. 
  \item Suppose $\chi_{lat}(G) = \chi(G \vee K_1) -1$ with a corresponding local antimagic total labeling $f$ of $G$. If $\sum_{i=1}^{p} f(v_i)\ne  w_f(v_j)$, $1\le j\le p$, then $\chi_{la}(G \vee K_1)=\chi(G \vee K_1)$. \end{enumerate} \end{theorem}

\ms\nt For $m\ge 2$ and $1\le i\le m$, let $G_i$ be a simple graph with an induced subgraph $H$. An {\it amalgamation} of $G_1,\ldots, G_m$ over $H$ is the simple graph obtained by identifying the vertices of $H$ of each $G_i$ so that the new obtained graph contains a subgraph $H$ induced by the identified vertices. Suppose $G$ is a graph with a proper subgraph $K_r$, $r\ge 1$. Let $A(mG, K_r)$ be the amalgamation of $m\ge 2$ copies of $G$ over $K_r$. Note that there may be many non-isomorphic $A(mG, K_r)$ graphs. For example, $A(2P_3,K_2)$ may be either $K_{1,3}$ or $P_4$. When $r=1$, the graph is also known as one-point union of graphs. Note that $A(mK_2,K_1)\cong K_{1,m}$ and $A(mK_3, K_1)$ is the friendship graph $f_m, m\ge 2$. In~\cite[Theorem 2.4]{Lau+SN}, the authors completely determined $\chi_{la}(A(mC_n, K_1))$ for $m\ge 2, n\ge 3$ where $2\le \chi(A(mC_n,K_1))\le 3$. Motivated by this, in this paper, we determine the $\chi_{lat}(A(mK_n,K_r))$ and $\chi_{la}(A(mK_n,K_r) \vee K_1)$ where $\chi(A(mK_n,K_r))=n$ for $m\ge1, n\ge 2, r\ge 0$. Sharp upper bounds on $\chi_{lat}(mK_n)$ are also obtained for odd $n\ge 3$.

\section{Amalgamations of Complete Graphs}\label{sec-chilat}

\nt Let $f$ be a total labeling of a simple $(p,q)$-graph $G$. Let $V(G)=\{u_1, \dots, u_p\}$. We define a total labeling matrix which is similar to the labeling matrix of an edge labeling introduced in \cite{Shiu+Lam+Lee2002}.

\ms\nt Suppose $f:V(G)\cup E(G) \rightarrow S$ is a mapping, where $S$ is a set of labels. A {\it total labeling matrix} $M$ of $f$ for $G$ is a $p\times p$ symmetric matrix in which the $(i,i)$-entry of $M$ is $f(u_i)$; the $(i,j)$-entry of $M$ is $f(u_iu_j)$ if $u_iu_j \in E$ and is $*$ otherwise. If $f$ is a local antimagic total labeling of $G$, then a total labeling matrix of $f$ is called a {\it local antimagic total labeling matrix} of $G$. Clearly the $i$-th row sum (and $i$-th column sum) is $w_f(u_i)$, where $*$'s are treated as zero. Thus the condition of a total labeling matrix being a local antimagic total labeling matrix is the $i$-th row sum different from the $j$-th row sum when $u_iu_j\in E$.

\ms\nt For $m\ge2$ and $n > r\ge 1$, let $V(A(mK_n,K_r))=\{v_{i,j}\,|\,1\le i\le m, 1\le j\le n\}$ and $E(A(mK_n, K_r)) = \{v_{i,j}v_{i,k}\,|\, 1\le i\le m, 1\le j < k\le n\}$, where $v_{1,j}=\cdots = v_{m,j}$ for each $n-r+1\le j\le n$.  For convenience, let $u_j=v_{1,j}$ for $n-r+1\le j\le n$. Note that $A(mK_n, K_r)\equiv mK_{n-r}\vee K_r$. We first list the vertices of the $m$ copies of $K_{n-r}$ in lexicographic order followed by $u_{n-r+1}, \dots, u_n$. Now let us show the structure of a total labeling matrix $M$ of the graph $A(mK_n, K_r)$ under this list of vertices as a block matrix. Namely,
\begin{equation}\label{eq-TLM} M=\begin{pmatrix}
 L_1 & \bigstar & \bigstar & \cdots & \cdots & \bigstar  & B_1\\
\bigstar & L_2 &\bigstar &\ddots & \cdots  &\bigstar & B_2\\
\vdots  &\ddots  & \ddots & \ddots &\ddots & \cdots &\vdots\\
\bigstar &\cdots  & \bigstar & L_i & \ddots&\bigstar  & B_i\\
\vdots  &\ddots  & \cdots & \ddots & \ddots & \vdots  &\vdots\\
\bigstar & \bigstar &\bigstar & \cdots & \bigstar & L_m & B_{m}\\
B_1^T & B_2^T & \cdots & B_i^T & \cdots & B_m^T & A\end{pmatrix},\end{equation}
where $L_i$ is an $(n-r)\times (n-r)$ symmetric matrix, $B_i$ is an $(n-r)\times r$ matrix, $1\le i\le m$, and $A$ is an $r\times r$ symmetric matrix. Here $\bigstar$ denotes an $(n-r)\times (n-r)$ matrix whose entries are $*$'s. Thus, the corresponding total labeling matrix of the $i$-th $K_n$ is
\begin{equation}\label{eq-TLMi} M_i=\begin{pmatrix}L_i & B_i\\
B_i^T & A\end{pmatrix}.\end{equation}

\nt Now, a local antimagic total labeling for the graph $A(mK_n, K_r)$ is obtained if we use the integers in $[1, N]$ for all entries of the upper triangular part of all $L_i$'s and $A$, and all entries of all $B_i$'s such that the row sums of each matrix $M_i$ are distinct, where $N=\frac{mn(n+1)}{2}-\frac{(m-1)r(r+1)}{2}$. We may extend the case to $r=0$. We let $A(mK_n, K_0)=mK_n$ by convention. For this case, all of $B_i$'s and $A$ in \eqref{eq-TLM} and \eqref{eq-TLMi} do not exist.

\ms\nt For a given matrix $B$, we shall use $\rs_i(B)$ and $\cs_j(B)$ to denote the $i$-th row sum and the $j$-th column sum of $B$, respectively. Also we shall use $\d(B)$ to denote the sum of the main diagonal of $B$ if $B$ is a square matrix. Suppose $S$ is a finite subset of $\Z$. Let $S^-$ and $S^+$ be a decreasing sequence and an increasing sequence of $S$, respectively.

\ms\nt We shall keep the notation defined above in this section.

\begin{lemma}\label{lem-rowsum}
Suppose $m\le n$. Let $M=(m_{j,k})$ be an $m\times n$ matrix with the following properties:
\begin{enumerate}[(a)]
\item $m_{k,j}=m_{j,k}$ for all $j,k$, where $1\le j< k\le m$;
\item $m_{j,j}<m_{k,k}$ if $1\le j<k\le m$;
\item for $j_1< k_1$ and $j_2< k_2$,  $(j_1,k_1)<(j_2, k_2)$ in lexicographic order implies that $m_{j_1,k_1}<m_{j_2, k_2}$.
\end{enumerate}
Then $\rs_j(M)$ is a strictly increasing function of $j$.
\end{lemma}

\begin{proof}\hspace*{\fill}{}
\vspace*{-7.5mm}
\begin{align*}
& \quad \rs_{j+1}(M)-\rs_j(M)\\&=\sum_{k=1}^{n} (m_{j+1,k}-m_{j,k})\\ & =\sum_{k=1}^{j-1} (m_{j+1,k}-m_{j,k}) + (m_{j+1,j}-m_{j,j}+m_{j+1,j+1}-m_{j,j+1})+\sum_{k=j+2}^{n}(m_{j+1,k}-m_{j,k})\\
& =\sum_{k=1}^{j-1} (m_{k,j+1}-m_{k,j}) + (m_{j+1,j+1}-m_{j,j})+\sum_{k=j+2}^{n}(m_{j+1,k}-m_{j,k}) >0.
\end{align*}
Note that the empty sum is treated as 0. This completes the proof.
\end{proof}

\begin{lemma}\label{lem-S}
For positive integers $t$ and $m$, let $S(a)=[m(a-1)+1, ma]$, $1\le a\le t$. We have:
\begin{enumerate}[(i)]
\item $\{S(a)\;|\; 1\le a\le t\}$ is a partition of $[1, mt]$.
\item If $a<b$, then every term of $S(a)$ is less than every term of $S(b)$.
\item For any $1\le a, b\le t$, the sum of the $i$-th term of $S^+(a)$ and that of $S^-(b)$ is independent of the choice of $i$, $1\le i\le m$.
\item For any $1\le a_l, b_l\le t$, $\sum\limits_{l=1}^k \left(i\mbox{-th term of }S^+(a_l)\right) + \sum\limits_{l=1}^k \left(i\mbox{-th term of } S^-(b_l)\right)$ is independent of the choice of $i$, $1\le i\le m$.
\end{enumerate}
\end{lemma}

\begin{proof}
\nt The first two parts are obvious. For (iii), the $i$-th terms of $S^+(a)$ and $S^-(b)$ are $m(a-1)+i$ and $m(b-1)+(m+1-i)$, respectively. So the sum is $m(a+b-1)+1$ which is independent of $i$. The last part follows from (iii).
\end{proof}

 \nt Before providing results about $\chi_{lat}(A(mK_n, K_r))$ for some $m,n,r$, we define a `sign matrix' $\S_n$ for even $n$.

\nt Let $\S_2=
\begin{pmatrix}
+1 & -1\\-1 & +1
\end{pmatrix}$ be a $2\times 2$ matrix and
$\S_4=\begin{pmatrix}\S_2 & \S_2\\
\S_2 & -\S_2
\end{pmatrix}$ be a $4\times 4$ matrix.
 Let $\S_{4k}$ be a $(4k)\times(4k)$ matrix given by the following block matrix, where $k\ge 2$:
\[\S_{4k}=\begin{pmatrix}
\S_4 & \cdots & \S_4\\
\vdots &\ddots & \vdots\\
\S_4 & \cdots & \S_4\end{pmatrix}.\]
\nt Let $\S_{4k+2}$ be a $(4k+2)\times(4k+2)$ matrix as the following block matrix, where $k\ge 1$:
\[\S_{4k+2}=\left(\begin{array}{c|l}
\S_{4k} & \begin{matrix}\S_2\\\vdots\\\S_2\end{matrix}\\\hline
\begin{matrix}\S_2 & \cdots & \S_2\end{matrix} &
\S_2\end{array}\right).\]

\nt We shall keep these notation in this section.

\begin{rem}\label{rem-sign-matrix} It is easy to see that each row and column sum of $\S_n$ are zero. Moreover, the diagonal sum of $\S_{4k}$ is zero.
\end{rem}


\begin{theorem}\label{thm-A(mKnKr)-even} For $m\ge 2$, $n$ even and $n>r\ge 0$,
$\chi_{lat}(A(mK_n, K_r))= n$.
\end{theorem}

\begin{proof}
Let $\S$ be the $(n-r)\times n$ matrix obtained from $\S_n$ by removing the last $r$ rows of $\S_n$. First, define an $(n-r)\times n$ matrix $\M'$ in which $(\M')_{j,k}=(\M')_{k,j}$ for $1\le j<k\le n-r$.
Assign the increasing sequence $[1, (n-r)(n+r-1)/2]$ in lexicographic order to the upper part of the off-diagonal entries of $\M'$, denoted $(j,k)$ if in row $j$ and column $k$ for $1\le j< k\le n$. Next, assign $[(n-r)(n+r-1)/2+1, (n-r)(n+r-1)/2+(n-r)]$ to the entries of the main diagonal of $\M'$ in natural order.

\ms\nt Now, define an $(n-r)\times n$ `guide matrix' $\M$ whose $(j,k)$-th entry is $(\S)_{j,k}(\M')_{j,k}$,  $1\le j\le n-r$ and $1\le k\le n$.

\ms \nt{\bf Stage 1:} We shall assign labels to the upper triangular entries of $L_i$'s and all the entries of $B_i$'s.  Note that if $r=0$, all of $B_i$'s and $A$ do not exist.  There are $N_1=(n-r)(n+r-1)/2+(n-r)=\frac{(n+r+1)(n-r)}{2}$ entries needed to be filled for each $i$.

\nt Now we shall use labels in $[1, mN_1]$ to fill in the $m$ submatrices $\begin{pmatrix}L_i & B_i\end{pmatrix}$, $1\le i\le m$. We use the sequences $S(a)$ defined in Lemma~\ref{lem-S}, where $t=N_1$. The $(j,k)$-entry of $M_i$ is the $i$-th term of $S^+(a)$ or $S^-(a)$ if the corresponding $(j,k)$-entry of $\M$ is $+a$ or $-a$ respectively, where $1\le j\le n-r$.

\ms \nt By Lemma~\ref{lem-S} (iv), $\rs_j(M_i)$ are the same for all $i$, $1\le i\le m$. In other words, $w(v_{i, j})$ is a constant function for a fixed $j$, $1\le j\le n-r$.

\ms\nt{\bf Stage 2: } Note that when $r=0$, the total labeling matrix $M$ does not have the last row and column of block matrices so that we only need to perform Stage~1 above. Thus, we now assume $r\ne 0$. Also note that all integers in $[1,mN_1]$ are used up in Stage 1. Use the increasing sequence $[mN_1+1, mN_1+r(r-1)/2]$ in lexicographic order for
the off-diagonal entries of $A$.  Lastly, use $[mN_1+r(r-1)/2+1, N]$ in natural order for the diagonals of $A$. The lower triangular part duplicates the upper triangular part.

\ms\nt Consider the matrix $M_m$. Clearly it satisfies the conditions of Lemma~\ref{lem-rowsum}. Hence $\rs_j(M_m)$ is a strictly increasing function of $j$, $1\le j\le n$. Thus $w(v_{i,j})=w(v_{m, j})$ is a strictly increasing function of $j$, $1\le j\le n-r$, for $1\le i\le m$.

\ms\nt By the structure of $B_i$ and $A$, and by Lemma~\ref{lem-S}(i), we have \begin{align*}w(u_{n-r+j}) & =\rs_{j}(A)+\sum\limits_{i=1}^{m}\rs_{j}(B_i^T)
=\rs_{j}(A)+\sum\limits_{i=1}^{m}\cs_{j}(B_i)\\&<
\rs_{j+1}(A)+\sum\limits_{i=1}^{m}\cs_{j+1}(B_i)=w(u_{n-r+j+1}),\end{align*} for $1\le j\le r-1$.

\ms \nt Thus $\chi_{lat}(A(mK_n, K_r))=n$ since $\chi(A(mK_n, K_r))=n$.
\end{proof}

\begin{rem}\label{rem-LAT-LA} Suppose $f$ is a local antimagic total labeling of a graph $G$ and $M$ is the corresponding total labeling matrix. From the proof of Theorem~\ref{thm-K1VG}(b) we can see that, if $\d(M)$ does not equal to every row (also column) sum of $M$, then $f$ induces a local antimagic labeling of $G\vee K_1$. Moreover, $\chi(G\vee K_1)\le \chi_{la}(G\vee K_1)\le \chi_{lat}(G)+1$.
\end{rem}

\begin{corollary}\label{cor-A(mKnKr)-even} For $m\ge 2$, $n$ even and $n>r\ge 0$,
$\chi_{la}(A(mK_{n+1}, K_{r+1}))=\chi_{la}(A(mK_{n}, K_{r})\vee K_1)= n+1$.
\end{corollary}
\begin{proof} Note that $A(mK_{n+1}, K_{r+1})\cong A(mK_{n}, K_{r})\vee K_1$. Since $\chi_{la}(A(mK_{n}, K_{r})\vee K_1)\ge \chi(A(mK_{n}, K_{r})\vee K_1)=n+1$, we only need to show $\chi_{la}(A(mK_{n}, K_{r})\vee K_1)\le n+1$. Keep the total labeling matrix $M$ of $A(mK_{n}, K_{r})$ in the proof of Theorem~\ref{thm-A(mKnKr)-even}.

\ms\nt Since each diagonal of $M$ is the largest entry in the corresponding column, $\d(M)$ is larger than each row sum of $M$. Thus, $\d(M)$ is greater than all vertex weights of $A(mK_{n}, K_{r})$. By Remark~\ref{rem-LAT-LA}, we get that $\chi_{la}(A(mK_{n}, K_{r})\vee K_1)\le n+1$.
\end{proof}

\begin{example}\label{ex-3K6} We take $m=3$, $n=6$ and $r=0$.

\begin{minipage}[b]{7.5cm}
\nt The guide matrix is\\
\fontsize{9}{10}\selectfont
$\M=\left(\begin{array}{*{6}{r}}
+16 & -1 & +2 & -3 & +4 & -5 \\
-1 & +17 & -6 & +7 & -8 & +9\\
+2 & -6 & -18 & +10 & +11 & -12\\
-3 & +7 & +10 & -19 & -13 & +14\\
+4 & -8 & +11 & -13 & +20 & -15\\
-5 & +9 & -12 & +14 & -15 & +21
\end{array}\right)$
\end{minipage}
\begin{minipage}[b]{8cm}
\fontsize{9}{10}\selectfont
$\begin{array}{c||*{6}{c|}|c}
M_1 & v_{1,1} & v_{1,2} & v_{1,3} & v_{1,4} & v_{1,5} & v_{1,6} &  \mbox{sum}\\\hline\hline
v_{1,1} & 46 & 3 & 4 & 9 & 10 & 15 & 87  \\\hline
v_{1,2} & 3 & 49 & 18 & 19 & 24 & 25 & 138 \\\hline
v_{1,3} & 4 & 18 & 54 & 28 & 31 & 36 & 171  \\\hline
v_{1,4} & 9 & 19 & 28 & 57 & 39 & 40 & 192  \\\hline
v_{1,5} & 10 & 24 & 31 & 39 & 58 & 45 & 207  \\\hline
v_{1,6} & 15 & 25 & 36 & 40 & 45 & 61 & 222 \\\hline
\end{array}$
\end{minipage}\\[3mm]

\nt\begin{minipage}[b]{8cm}
\fontsize{9}{10}\selectfont
$\begin{array}{c||*{6}{c|}|c}
M_2 & v_{2,1} & v_{2,2} & v_{2,3} & v_{2,4} & v_{2,5} & v_{2,6} & \mbox{sum}\\\hline\hline
v_{2,1} & 47 & 2 & 5 & 8 & 11 & 14 & 87  \\\hline
v_{2,2} & 2 & 50 & 17 & 20 & 23 & 26 & 138 \\\hline
v_{2,3} & 5 & 17 & 53 & 29 & 32 & 35 & 171  \\\hline
v_{2,4} & 8 & 20 & 29 & 56 & 38 & 41 & 192  \\\hline
v_{2,5} & 11 & 23 & 32 & 38 & 59 & 44 & 207  \\\hline
v_{2,6} & 14 & 26 & 38 & 41 & 44 & 62 & 222 \\\hline
\end{array}$
\end{minipage}
\begin{minipage}[b]{8cm}
\fontsize{9}{10}\selectfont
$\begin{array}{c||*{6}{c|}|c}
M_3 & v_{3,1} & v_{3,2} & v_{3,3} & v_{3,4} & v_{3,5} & v_{3,6} & \mbox{sum}\\\hline\hline
v_{3,1} & 48 & 1 & 6 & 7 & 12 & 13 & 87  \\\hline
v_{3,2} & 1 & 51 & 16 & 21 & 22 & 27 & 138 \\\hline
v_{3,3} & 6 & 16 & 52 & 30 & 33 & 34 & 171  \\\hline
v_{3,4} & 7 & 21 & 30 & 55 & 37 & 42 & 192  \\\hline
v_{3,5} & 12 & 22 & 33 & 37 & 60 & 43 & 207  \\\hline
v_{3,6} & 13 & 27 & 34 & 42 & 43 & 63 & 222 \\\hline
\end{array}$
\end{minipage}

\ms\nt The above matrices give $\chi_{lat}(3K_6) = 6$. Let $v$ be the vertex of $K_1$. If the main diagonal labels are the edge labels of $3K_6 \vee K_1$ incident with $v$, then the induced label of $v$ is $981$. Thus, the matrices give $\chi_{la}(3K_6 \vee K_1)=7$.

\ms\nt Deleting edge $vv_{3,6}$ of label 63 from $3K_6\vee K_1$, we get a local antimagic labeling of $(3K_6\vee K_1)-vv_{3,6}$. Thus, by symmetry, $\chi_{la}((3K_6\vee K_1)-e)=7$ for $e$ not belonging to any $K_6$.

\ms \nt Delete the edge $v_{3,1}v_{3,2}$ that has label 1 and reduce all other labels by 1. We get $\chi_{lat}(3K_6 - e)=6$ by symmetry for $e$ that belongs to any $K_6$. Now, if the main diagonal labels are the edge labels of $(3K_6-e) \vee K_1$ incident with $v$, then we have $\chi_{la}(3K_6\vee K_1) - e)=7$ for $e$ that belongs to any $K_6$. \rsq

\end{example}

\begin{example} We take $m=3$, $n=6$ and $r=1$. The guide matrix is obtained from the guide matrix of Example~\ref{ex-3K6} by deleting the last row. So we have
\begin{align*}
M_1& =\left(\begin{array}{*{5}{c}|c||c}
46 & 3 & 4 & 9 & 10 & 15 & 87  \\
3 & 49 & 18 & 19 & 24 & 25 & 138 \\
4 & 18 & 54 & 28 & 31 & 36 & 171  \\
9 & 19 & 28 & 57 & 39 & 40 & 192  \\
10 & 24 & 31 & 39 & 58 & 45 & 207  \\\hline
15 & 25 & 36 & 40 & 45 & 61 & 222
\end{array}\right),\quad
M_2=\left(\begin{array}{*{5}{c}|c||c}
47 & 2 & 5 & 8 & 11 & 14 & 87  \\
2 & 50 & 17 & 20 & 23 & 26 & 138 \\
5 & 17 & 53 & 29 & 32 & 35 & 171  \\
8 & 20 & 29 & 56 & 38 & 41 & 192  \\
11 & 23 & 32 & 38 & 59 & 44 & 207  \\\hline
14 & 26 & 38 & 41 & 44 & 61 & 221
\end{array}\right),\\
M_3& =
\left(\begin{array}{*{5}{c}|c||c}
48 & 1 & 6 & 7 & 12 & 13 & 87  \\
1 & 51 & 16 & 21 & 22 & 27 & 138 \\
6 & 16 & 52 & 30 & 33 & 34 & 171  \\
7 & 21 & 30 & 55 & 37 & 42 & 192  \\
12 & 22 & 33 & 37 & 60 & 43 & 207  \\\hline
13 & 27 & 34 & 42 & 43 & 61 & 220
\end{array}\right).
\end{align*}
The last column of each matrix is the corresponding row sum. Now $w(u_6)=222+221+220-2\times 61=541$. Hence $\chi_{lat}(A(3K_6, K_1))=6$.
Since $\d(M)=856$, $\chi_{la}(A(3K_6, K_1)\vee K_1)=7$.
\rsq
\end{example}

\begin{corollary}\label{cor-LA-A(mK4k)} For $m\ge 2$,
$\chi_{la}(mK_{4k+1})= 4k+1$.
\end{corollary}
\begin{proof}
Consider the total labeling matrix of $mK_{4k}$ defined in the proof of Theorem~\ref{thm-A(mKnKr)-even}. For each matrix $M_i=L_i$, $1\le i\le m$, we add the $(n+1)$-st extra column at the right of $M_i$ with entry $*$. For each row of this matrix, swap the diagonal entry with the entry of the $(n+1)$-st column. Add the $(n+1)$-st extra row to this matrix and let the $(n+1, n+1)$-entry be $*$ and then make the resulting matrix $Q_i$ to be symmetric. Then $Q_i$ is a labeling matrix of the $i$-th copy of $K_{4k+1}$.

\ms\nt By Remark~\ref{rem-sign-matrix} and Lemma~\ref{lem-S} (iv), all the diagonal sums of $M_i$'s are the same, $1\le i\le n$. Thus the $j$-th row sum of $Q_i$ is independent of $i$, $1\le j\le n+1$. Hence we have $\chi_{la}(mK_{4k+1})\le 4k+1$. Since $\chi(mK_{4k+1})=4k+1$, $\chi_{la}(mK_{4k+1})=4k+1$.
\end{proof}

\begin{example} We take $m=3$, $n=4$. So
\[\M =\begin{pmatrix}
+7 & -1 & +2 & -3\\
-1 & +8 & -4 & +5\\
+2 & -4 & -9 & +6\\
-3 & +5 & +6 & -10\end{pmatrix},\]
and
\[M_1=\left(\begin{array}{cccc||c}
19 & 3 & 4 & 9 & 35\\
3 & 22 & 12 & 13 & 50\\
4 & 12 & 27 & 16 & 59\\
9 & 13 & 16 & 30 & 68\end{array}\right),
M_2=\left(\begin{array}{cccc||c}
20 & 2 & 5 & 8 & 35\\
2 & 23 & 11 & 14 & 50\\
5 & 11 & 26 & 17 & 59\\
8 & 14 & 17 & 29 & 68\end{array}\right),
M_3=\left(\begin{array}{cccc||c}
21 & 1 & 6 & 7 & 35\\
1 & 24 & 10 & 15 & 50\\
6 & 10 & 25 & 18 & 59\\
7 & 15 & 18 & 28 & 68\end{array}\right)
\]
Thus $\chi_{lat}(3K_4)=4$.

\nt Let
\begin{align*} Q_1& =\left(\begin{array}{ccccc||c}
* & 3 & 4 & 9 & 19 & 35\\
3 & * & 12 & 13 & 22 & 50\\
4 & 12 & * & 16 &  27 & 59\\
9 & 13 & 16 & * & 30 & 68\\
19 & 22 & 27 & 30 & * & 98\end{array}\right),
Q_2=\left(\begin{array}{ccccc||c}
* & 2 & 5 & 8 & 20 & 35\\
2 & * & 11 & 14 & 23 & 50\\
5 & 11 & * & 17 & 26 & 59\\
8 & 14 & 17 & * & 29 & 68\\
20 & 23 & 26 & 29 & * & 98\end{array}\right),\\
Q_3 & =\left(\begin{array}{ccccc||c}
* & 1 & 6 & 7 & 21 & 35\\
1 & * & 10 & 15 & 24 & 50\\
6 & 10 & * & 18 & 25 & 59 \\
7 & 15 & 18 & * & 28 & 68\\
21 & 24 & 25 & 28 & * & 98
\end{array}\right).
\end{align*}
Thus $\chi_{la}(3K_5)=5$.
\rsq\end{example}

\begin{theorem}\label{thm-A(mKnKr)-odd} For $m\ge 2$, odd $n\ge 3$ and $n>r\ge 3$,
$\chi_{lat}(A(mK_n, K_r))= n$.
\end{theorem}

\begin{proof}
Suppose $r$ is odd so that $n-r$ is even.

\nt{\bf Stage 1: } Using the same approach of the proof of Theorem~\ref{thm-A(mKnKr)-even}, we construct an $(n-r)\times (n-r)$ guide matrix $\M$.

\ms\nt Similar to the proof of Theorem~\ref{thm-A(mKnKr)-even}, we use the guide matrix $\M$ for all entries of $L_i$, $1\le i\le m$, using labels in $[1, mN_2]$, where $N_2=(n-r)(n-r+1)/2$. Thus, $\rs_j(L_i)$ is a function only depending on $j$ and is strictly increasing, for $1\le i\le m$ and $1\le j\le n-r$.

\ms \nt{\bf Stage 2: }
Use $[mN_2+1, mN_2+(m(n-r)+1)r]$ to form an $(m(n-r)+1)\times r$ magic rectangle $\Omega$. Note that the existence of this magic rectangle is referred to in \cite{Hagedorn}. Let
\[\begin{pmatrix} B_1\\ B_2\\\vdots\\B_m\\\alpha\end{pmatrix}=\Omega,\] where $\alpha=(A_{1,1}, A_{2,2}, \dots, A_{r,r})$. Now, for a fixed $i$, $w(v_{i,j})=\rs_j(L_i)+\rs_j(B_i)$. So $w(v_{i,j})$ is a function only depending on $j$ and is strictly increasing, for $1\le i\le m$ and $1\le j\le n-r$.

\ms \nt{\bf Stage 3: }
Use the increasing sequence $[mN_2+(m(n-r)+1)r+1, N]$ in lexicographic order for the remaining entries of the upper triangular part of $A$. For $1\le j\le r$, \begin{align*}w(u_{n-r+j}) & = \sum_{i=1}^{m} \rs_j(B^T_i) +
\rs_j(A) =\sum_{i=1}^{m} \cs_j(B_i)+A_{j,j} +\sum_{\begin{smallmatrix}l=1\\l\ne j\end{smallmatrix}}^{r}A_{j,l}\\ & =
\cs_j(\Omega) + \sum_{l=1}^{j-1}A_{j,l} +A_{j, j+1}+\sum_{l=j+2}^{r}A_{j,l}
\\ & <
\cs_{j}(\Omega) + \sum_{l=1}^{j-1}A_{j+1,l} +A_{j+1, j}+\sum_{l=j+2}^{r}A_{j+1,l}\tag{since $r\ge 3$, there is at least one non-empty sum}\\
& = \cs_{j+1}(\Omega) + \sum_{\begin{smallmatrix}l=1\\l\ne j+1\end{smallmatrix}}^{r}A_{j+1,l} = w(u_{n-r+j+1}). \end{align*}
So $w(u_{n-r+j})$ is a strictly increasing function of $j$ for $1\le j\le r$.

\begin{align*} w(v_{m, n-r})&= \rs_{n-r}(L_m)+\rs_{n-r}(B_m) =\cs_{n-r}(L_m)+\rs_1(\alpha) \tag{since $\Omega$ is a magic rectangle}\\
& <\cs_1(B_m) +\rs_1(\alpha) < \rs_1(B_m^T) +\sum_{k=1}^{r}A_{1,k} <w(u_{n-r+1}).
\end{align*}

\ms\nt Thus, we have $\chi_{lat}A(mK_n, K_r)=n$.

\ms\nt Suppose $r$ is even so that $n-r+1$ is even.

\begin{enumerate}[(a)]

\item Suppose $m$ is odd so that $m(n-r)$ and $r-1$ are odd and at least 3.

\nt{\bf Stage 1(a): }  We use a  modification of the proof of Theorem~\ref{thm-A(mKnKr)-even}. Firstly we use the $(n-r)\times (n-r+1)$ sign matrix $\S$. Next, we define an $(n-r)\times (n-r+1)$ matrix $\M'$ by assigning the increasing sequence $[1, (n-r)(n-r+1)/2]$ in lexicographic order to the off-diagonal entries of the upper triangular part of $\M'$. Now, assign $[(n-r)(n+r+1)/2-(n-r-1), (n-r)(n+r+1)/2]$ to the entries of the main diagonal of $\M'$ in natural order. The $(n-r)\times (n-r+1)$ guide matrix $\M$ is defined the same way as in the proof of Theorem~\ref{thm-A(mKnKr)-even}.

\ms\nt Write $B_i=\begin{pmatrix}\beta_i & X_i\end{pmatrix}$, where $\beta_i$ and $X_i$ are $(n-r)\times 1$ and $(n-r)\times (r-1)$ matrices, respectively. Similar to Stage~1 of the odd $r$ case, for $\begin{pmatrix}L_i & \beta_i\end{pmatrix}$ use the labels in $[1, mZ_1]\cup [m[Z_2-(n-r)]+1, mZ_2]$, $1\le i\le m$, where $Z_1=(n-r)(n-r+1)/2$ and $Z_2=Z_1+(n-r)r$.

\nt Thus, $\rs_j\begin{pmatrix}L_i & \beta_i\end{pmatrix}$ is a function only depending on $j$ and is a strictly increasing function of $j$, $1\le i\le m$ and $1\le j\le n-r$.

\ms \nt{\bf Stage 2(a): } Use $[mZ_1+1, m[Z_1+(n-r)(r-1)]]$ to form an $m(n-r)\times (r-1)$ magic rectangle $\Omega$.
Let
\[\begin{pmatrix} X_1\\ X_2\\\vdots\\X_m\end{pmatrix}=\Omega.\]

\nt Now, for a fixed $i$, $w(v_{i,j})=\rs_j\begin{pmatrix}L_i & \beta_i\end{pmatrix}+\rs_j(X_i)$. So $w(v_{i,j})$ is a function only depending on $j$ and is a strictly increasing function of $j$ for $1\le j\le n-r$ and $1\le i\le m$.

\ms \nt{\bf Stage 3(a): }  Use the increasing sequence $[mZ_2+1, mZ_2+(r-1)r/2]$ in lexicographic order for the off-diagonal entries of the upper triangular part of $A$ and then use $[mZ_2+(r-1)r/2+1, N]$ in natural order for the diagonals of $A$. By Lemma~\ref{lem-rowsum}, $\rs_j(A)$ is a strictly increasing function of $j$ for $1\le j\le r$.

\nt Now, for $2\le j\le r$, \begin{align*}w(u_{n-r+j}) & = \sum_{i=1}^{m}\sum_{k=1}^{n-r}(B_i^T)_{j, k} +\sum_{l=1}^{r}A_{j,l} = \sum_{i=1}^{m}\sum_{k=1}^{n-r}(B_i)_{k, j} +\sum_{l=1}^{r}A_{j,l}\\ & = \cs_{j-1}(\Omega)+\rs_j(A). \end{align*}

\nt So $w(u_{n-r+j})$ is a strictly increasing function of $j$ for $2\le j\le r$.

\nt Next
\begin{align*}w(u_{n-r+1}) & = \sum_{i=1}^{m}\sum_{k=1}^{n-r}(B_i^T)_{1, k} +\sum_{l=1}^{r}A_{1,l} = \sum_{i=1}^{m}\sum_{k=1}^{n-r}(B_i)_{k, 1} +\rs_1(A)\\
& <\sum_{i=1}^{m}\sum_{k=1}^{n-r}(B_i)_{k, 2} +\rs_2(A)=w(u_{n-r+2}). \end{align*}

Now
\begin{align*} w(v_{m, n-r})&= \rs_{n-r}(L_m)+\rs_{n-r}(B_m) =\cs_{n-r}(L_m)+ \rs_{n-r}(B_m)\\
& <\cs_1(B_m) +\sum_{k=1}^{r}A_{1,k} =\rs_1(B_m^T) +\sum_{k=1}^{r}A_{1,k} <w(u_{n-r+1}).
\end{align*}

\nt Thus we have $\chi_{lat}(A(mK_n, K_r))=n$.

\item Suppose $m$ is even so that $m(n-r)$ is even.

\nt{\bf Stage 1(b): } Similar to Stage~1 of the odd $r$ case we define an $(n-r)\times (n-r+1)$ guide matrix $\M$.

\ms\nt Write $B_i=\begin{pmatrix}\beta_i & X_i\end{pmatrix}$, where $\beta_i$ and $X_i$ are $(n-r)\times 1$ and $(n-r)\times (r-1)$ matrices, respectively. Similar to Stage~1 of the odd $r$ case, for $\begin{pmatrix}L_i & \beta_i\end{pmatrix}$ use labels in $[1, mN_3]$, $1\le i\le m$, where $N_3=(n-r+3)(n-r)/2$. Thus, $\rs_j\begin{pmatrix}L_i & \beta_i\end{pmatrix}$ is a function only depending on $j$ and is a strictly increasing function of $j$, $1\le i\le m$ and $1\le j\le n-r$.

\ms \nt{\bf Stage 2(b): } $[mN_3+2, mN_3+(m(n-r)+1)(r-1)+1]$ to form an $(m(n-r)+1)\times (r-1)$ magic rectangle $\Omega$. We will assign $mN_3+1$ to $A_{1,1}$ in the next stage.


Let
\[\begin{pmatrix} X_1\\ X_2\\\vdots\\X_m\\\alpha\end{pmatrix}=\Omega,\]
where $\alpha$ is a $1\times (r-1)$ matrix.

\nt Same as Stage~2(a), we have $w(v_{i,j})$ is a function only depending on $j$ and is a strictly increasing function for $1\le j\le n-r$ and $1\le i\le m$.

\ms \nt{\bf Stage 3(b): } Let $A_{1,1}=mN_3+1$. Use the increasing sequence $[mN_3+(m(n-r)+1)(r-1)+2, N]$  in lexicographic order for the remaining entries of the upper triangular part of $A$. 

By the same proof of Stage~3 of the odd $r$ case, we have $w(u_{n-r+j})$ is a strictly increasing function of $j$, $2\le j\le r$. By a similar proof of Stage~3(a) we have $w(u_{n-r+1})<w(u_{n-r+2})$.

Now
\begin{align} w(v_{m, n-r})&= \rs_{n-r}(L_m)+\rs_{n-r}(B_m) =
\sum_{k=1}^{n-r} (L_m)_{n-r,k}+ (\beta_{m})_{n-r,1} + \sum_{k=1}^{r-1}(X_m)_{n-r, k}\nonumber\\
&= \sum_{k=1}^{n-r-1} (L_m)_{k,n-r}+ (L_m)_{n-r,n-r} +(\beta_{m})_{n-r,1} + \sum_{k=1}^{r-1}(X_m)_{n-r, k}\label{eq-odd-n}\\
& <\sum_{k=1}^{n-r-1} (\beta_m)_{k,1}+ A_{1,1} +(\beta_{m})_{n-r,1} + \sum_{k=2}^{r}A_{1, k}<w(u_{n-r+1}).\nonumber
\end{align}

\nt Thus we have $\chi_{lat}(A(mK_n, K_r))=n$.
\end{enumerate}
\end{proof}

\begin{corollary}For $m\ge 2$, $n$ odd and $n>r\ge 3$,
$\chi_{la}(A(mK_{n+1}, K_{r+1}))=\chi_{la}(A(mK_{n}, K_{r})\vee K_1)= n+1$.
\end{corollary}

\begin{proof} Similar to Corollary~\ref{cor-A(mKnKr)-even}, we only need to show that $\chi_{la}(A(mK_{n}, K_{r})\vee K_1)\le n+1$. Keep the construction in the proof of Theorem~\ref{thm-A(mKnKr)-odd}.

\begin{enumerate}[(A)]
\item  Suppose $r$ is odd so that $n-r \ge 2$ is even.
\begin{align*}w(v_{m, n-r}) & = \sum_{i=1}^{n-r} (L_{m})_{n-r,i} +\sum_{i=1}^{r} (B_{m})_{n-r,i} \\ &=\sum_{i=1}^{n-r} (L_{m})_{i,n-r} +\sum_{i=1}^r A_{i,i}\tag{since $\Omega$ is a magic rectangle}\\
&
<\sum_{i=1}^{n-r} (L_m)_{i,i} +\sum_{i=1}^r A_{i,i}  \le \d(M).
\end{align*}
\nt Next
\begin{align*}
\d(M)=&  \sum_{i=1}^{m}\sum_{j=1}^{n-r}(L_i)_{j,j} + A_{1,1}+\sum_{k=2}^{r}A_{k,k}
< \sum_{i=1}^{m}\sum_{j=1}^{n-r}(B_i)_{j,1} + A_{1,1}+\sum_{k=2}^{r} A_{k,k}\\
& < \sum_{i=1}^{m}\sum_{j=1}^{n-r}(B_i)_{j,1} + A_{1,1}+\sum_{k=2}^{r} A_{1,k}
\\&= w(u_{n-r+1}).
\end{align*}
Thus we have $w(v_{m, n-r})<\d(M)< w(u_{n-r+1})$.
\item Suppose $r$ is even and $m$ is odd. Since each diagonal of $M$ is the largest entry in the corresponding column, $\d(M)$ is larger than all the vertex weights.
\item Suppose $r$ is even and $m$ is even.
\begin{align*}
\d(M) & =  \sum_{i=1}^{m}\sum_{j=1}^{n-r}(L_i)_{j,j} + \sum_{k=1}^{r}A_{k,k}
< \sum_{i=1}^{m}\sum_{j=1}^{n-r}(X_i)_{j,1}+ \sum_{k=1}^{r}A_{2,k} \\
& =\sum_{i=1}^{m}\sum_{j=1}^{n-r}(X_i^T)_{1,j}+ \sum_{k=1}^{r}A_{2,k}=w(u_{n-r+2}).
\end{align*}
\begin{align*}
w(u_{n-r+1})&=\sum_{i=1}^{m}\sum_{k=1}^{n-r} (\beta_i)_{k,1}+ A_{1,1} + \sum_{k=2}^{r}A_{1, k}\\
&<\sum_{i=1}^{m}\sum_{k=1}^{n-r} (L_i)_{k,k}+A_{1,1} + \sum_{k=2}^{r}A_{1, k}\\
& <\sum_{i=1}^{m}\sum_{k=1}^{n-r} (L_i)_{k,k}+A_{1,1} + \sum_{k=2}^{r}A_{k, k}=\d(M)
\end{align*}

\end{enumerate}
\nt From Remark~\ref{rem-LAT-LA}, we have $\chi_{la}(A(mK_{n}, K_{r})\vee K_1)\le n+1$.
\end{proof}

\begin{example} We take $m=2$, $n=7$ and $r=3$. The guide matrix is
\[\M=\left(\begin{array}{cccc}
+7 & -1 & +2 & -3 \\
-1 & +8 & -4 & +5\\
+2 & -4 & -9 & +6\\
-3 & +9 & +6 & -10
\end{array}\right);\quad
L_1 =\left(\begin{array}{cccc}
13 & 2 & 3 & 6 \\
2 & 15 & 8 & 9\\
3 & 8 & 18 & 11\\
6 & 9 & 11 & 20
\end{array}\right),\quad
L_2=\left(\begin{array}{cccc}
14 & 1 & 4 & 5\\
1 & 16 & 7 & 10\\
4 & 7 & 17 & 12\\
5 & 10 & 12 & 19
\end{array}\right).
\]
\[\Omega=\begin{pmatrix} B_1\\B_2\\\alpha\end{pmatrix}=\left(\begin{array}{ccc}
25 & 30 & 47\\
34 & 39 & 29\\
36 & 44 & 22\\
41 & 28 & 33\\\hline
43 & 21 & 38\\
45 & 26 & 31\\
23 & 37 & 42\\
27 & 35 & 40\\\hline
32 & 46 & 24\end{array}
\right),\quad
A=\left(\begin{array}{ccc}
32 & 48 & 49\\
48 & 46 & 50\\
49 & 50 & 24\end{array}\right).
\]
\[M=\left(\begin{array}{*{4}{c}|*{4}{c}|*{3}{c}||c}
13 & 2 & 3 & 6 & * & * & * & * & 25 & 30 & 47 & 126\\
2 & 15 & 8 & 9 & * & * & * & * & 34 & 39 & 29 & 136\\
3 & 8 & 18 & 11 & * & * & * & * & 36 & 44 & 22 & 142\\
6 & 9 & 11 & 20 & * & * & * & * & 41 & 28 & 33 & 148\\\hline
* & * & * & * & 14 & 1 & 4 & 5 & 43 & 21 & 38 & 126\\
* & * & * & * & 1 & 16 & 7 & 10 & 45 & 26 & 31 & 136\\
* & * & * & * & 4 & 7 & 17 & 12 & 23 & 37 & 42 & 142\\
* & * & * & * & 5 & 10 & 12 & 19 & 27 & 35 & 40 & 148\\\hline
25 & 34 & 36 & 41 & 43 & 45 & 23 & 27 & 32  & 48 & 49 & 403\\
30 & 39 & 44 & 28 & 21 & 26 & 37 & 35 & 48 & 46  & 50 & 404\\
47 & 29 & 22 & 33 & 38 & 31 & 42 & 40 & 49 & 50 & 24  & 405
\end{array}\right).\]
\nt Thus $\chi_{lat}(A(2K_7, K_3))=7$. Since $\d(M)=234$, $\chi_{la}(A(2K_7, K_3)\vee K_1)=8$. \rsq
\end{example}

\begin{example}\label{ex-3K7-4} We take $m=3$, $n=7$ and $r=4$. The guide matrix is
\[\M=\left(\begin{array}{cccc}
+16 & -1 & +2 & -3 \\
-1 & +17 & -4 & +5\\
+2 & -4 & -18 & +6
\end{array}\right).\]
\[(L_1|\beta_1) =\left(\begin{array}{ccc|c}
46 & 3 & 4 & 9\\
3 &	49 & 12 & 13\\
4 & 12 & 54 & 16
\end{array}\right),\quad
(L_2|\beta_2)=\left(\begin{array}{ccc|c}
47 & 2 & 5 & 8\\
2 & 50 & 11 & 14\\
5 & 11 & 53 & 17
\end{array}\right),\quad
(L_3|\beta_3)=\left(\begin{array}{ccc|c}
48 & 1 & 6 & 7\\
1 & 51 & 10 & 15\\
6 & 10 & 52 & 18
\end{array}\right).\]

\[\Omega=\begin{pmatrix}X_1\\X_2\\X_3\end{pmatrix}=\left(\begin{array}{ccc}
27 & 32 & 37\\
20 & 34 & 42\\
31 & 39 & 26\\\hline
36 & 41 & 19\\
29 & 43 & 24\\
40 & 21 & 35\\\hline
45 & 23 & 28\\
38 & 25 & 31\\
22 & 30 & 44\end{array}
\right),\quad
A=\left(\begin{array}{cccc}
61 & 55 & 56 & 57\\
55 & 62 & 58 & 59\\
56 & 58 & 63 & 60\\
57 & 59 & 60 & 64\end{array}\right).\]
\[M=\left(\begin{array}{*{3}{c}|*{3}{c}|*{3}{c}|*{4}{c}||c}
46 & 3 & 4 & * & * & * & * & * & * & 9 & 27 & 32 & 37 & 158\\
3 &	49 & 12 & * & * & * & * & * & * & 13 & 20 & 34 & 42 & 173\\
4 & 12 & 54 & * & * & * & * & * & * & 16 & 31 & 39 & 26 & 182\\\hline
* & * & * & 47 & 2 & 5 & * & * & * & 8 & 36 & 41 & 19 & 158\\
* & * & * & 2 & 50 & 11 & * & * & * & 14 & 29 & 43 & 24 & 173\\
* & * & * & 5 & 11 & 53 & * & * & * & 17 & 40 & 21 & 35 & 182\\\hline
* & * & * & * & * & * & 48 & 1 & 6 & 7 & 45 & 23 & 28 & 158\\
* & * & * & * & * & * & 1 & 51 & 10 & 15 & 38 & 25 & 33 & 173\\
* & * & * & * & * & * & 6 & 10 & 52 & 18 & 22 & 30 & 44 & 182\\\hline
9 & 13 & 16 & 8 & 14 & 17 & 7 & 15 & 18 & 61 & 55 & 56 & 57 & 346\\
27 & 20 & 31 & 36 & 29 & 40 & 45 & 38 & 22 & 55 & 62 & 58 & 59 & 522\\
32 & 34 & 39 & 41 & 43 & 21 & 23 & 25 & 30 & 56 & 58 & 63 & 60 & 525\\
37 & 42 & 26 & 19 & 24 & 35 & 28 & 33 & 44 & 57 & 59 & 60 & 64 & 528
\end{array}\right).\]
Thus $\chi_{lat}(A(3K_7, K_4))=7$. Since $\d(M)=700$, $\chi_{la}(A(3K_7, K_4)\vee K_1)=8$. \rsq
\end{example}

\begin{example} We take $m=2$, $n=7$ and $r=4$. The guide matrix is
\[\M=\left(\begin{array}{cccc}
+7 & -1 & +2 & -3 \\
-1 & +8 & -4 & +5\\
+2 & -4 & -9 & +6
\end{array}\right),\quad
(L_1|\beta_1) =\left(\begin{array}{ccc|c}
13 & 2 & 3 & 6\\
2 &	15 & 8 & 9\\
3 & 8 & 18 & 11
\end{array}\right),\quad
(L_2|\beta_2)=\left(\begin{array}{ccc|c}
14 & 1 & 4 & 5\\
1 & 16 & 7 & 10\\
4 & 7 & 17 & 12
\end{array}\right).\]
\[\Omega=\begin{pmatrix}X_1\\X_2\\\alpha\end{pmatrix}=\left(\begin{array}{ccc}
27 & 25 & 38\\
24 & 31 & 35\\
21 & 32 & 37\\\hline
39 & 22 & 29\\
36 & 26 & 28\\
33 & 34 & 23\\\hline
30 & 40 & 20\end{array}
\right),\quad
A=\left(\begin{array}{cccc}
19 & 30 & 40 & 20\\
30 & 41 & 42 & 43 \\
40 & 42 & 44 & 45\\
20 & 43 & 45 & 46\end{array}\right).
\]
\[M=\left(\begin{array}{*{3}{c}|*{3}{c}|*{4}{c}||c}
13 & 2 & 3 & * & * & * & 6 & 27 & 25& 38  & 114\\
2 &	15 & 8 & * & * & * & 9 & 24 & 31 & 35 & 124\\
3 & 8 & 18 & * & * & * & 11 & 21 & 32 & 37 & 130\\\hline
* & * & * & 14 & 1 & 4 & 5 & 39 & 22 & 29 & 114\\
* & * & * & 1 & 16 & 7 & 10 & 36 & 26 & 28 & 124\\
* & * & * & 4 &	7 & 17 & 12 & 33 & 34 & 23 & 130\\\hline
6 & 9 & 11 & 5 & 10 & 12 & 19 & 30 & 40 & 20 & 162\\
27 & 24 & 21 & 39 & 36 & 33 & 30 & 41 & 42 & 43 & 336\\
25 & 31 & 32 & 22 & 26 & 34 & 40 & 42 & 44 & 45 & 341\\
38 & 35 & 37 & 29 & 28 & 23 & 20 & 43 & 45 & 46 & 344
\end{array}\right).\]
Thus $\chi_{lat}(A(2K_7, K_4))=7$. Since $\mathcal D(M)=243$, $\chi_{la}(A(2K_7, K_4)\vee K_1)=8$. \rsq
\end{example}

\begin{theorem}\label{thm-mKn-3}
For $m\ge 2$ and $n=4k+3\ge 3$, $\chi_{lat}(mK_n)\le n+1$.
\end{theorem}

\begin{proof}
Consider $\S_{4k+2}$. Change each diagonal entry of $\S_{4k+2}$ from $\pm 1$ to $\pm 2$. Let this new sign matrix be $\S'$. Now define $\S_{4k+3}$ by extending $\S'$ to a $4k+3$ square matrix by adding the last column and row. Use $+1$, $-1$ or $+2$ for each entry of this new column and row such that the row sum and the column sums of $\S_{4k+3}$ are zero.

\ms\nt Define a symmetric matrix $\M'$ of order $4k+3$ by using the increasing sequence $[1, (2k+1)(4k+3)]$ in lexicographic order for the upper triangular entries of $\M'$, and use $*$ for all diagonal entries.

\nt Now let $\M$ be the guide matrix whose $(j,l)$-th entry is $*$ if $j=l$; and $(\S_{4k+3})_{j,l}(\M')_{j,l}$, $1\le j,l\le 4k+3$ if $j\ne l$.

\ms \nt{\bf Stage 1:} Using the same procedure as Stage~1 in the proof of Theorem~\ref{thm-A(mKnKr)-even} fill the off-diagonals of the $m$ submatrices $L_i$, $1\le i\le m$, with labels in $[1, mN_4]$, where $N_4=\frac{n(n-1)}{2}$.

\ms \nt{\bf Stage 2:} Let $T(2a-1)=\{mN_4+2l-1+2m(a-1)\;|\; 1\le l\le m\}$ and $T(2a)=\{mN_4+2l+2m(a-1)\;|\; 1\le l\le m\}$,  $1\le a\le 2k+1$. For $1\le j\le 4k+2$, the $(j,j)$-entry of $L_i$ is the $i$-th term of $T^-(j)$ or $T^+(j)$, if the corresponding $(j,j)$-entry of $\mathcal M$ is $-2$ or $+2$, respectively.

\ms \nt{\bf Stage 3:}  Let $U(1)=\{mN_4+(4k+2)m+2l-1\;|\; 1\le l\le \lceil m/2\rceil\}$ and $U(2)=\{mN_4+(4k+2)m+2l\;|\; 1\le l\le \lfloor m/2\rfloor\}$. Let $U$ be the compound sequence $U^+(1) U^+(2)$, i.e., list the terms of $U^+(1)$ first and then follow with the terms of $U^+(2)$.
The $(4k+3,4k+3)$-entry of $L_i$  is the $i$-th term of $U$.

\ms\nt
For a fixed $j$, $1\le j\le 4k+2$, according to the structures of $S(a)$'s and $T(a)$'s, $(L_{i+1})_{j, l}-(L_i)_{j,l}$ is $(\S_{4k+3})_{j,l}$ for $1\le i\le m-1$, $1\le l\le 4k+3$. Thus, $w(v_{i,j})$ is a constant for a fixed $j$.

\nt Similarly, $(L_{i+1})_{4k+3, l}-(L_i)_{4k+3,l}$ is $(\S_{4k+3})_{4k+3,l}$ for $1\le i\le m-1$, $1\le l\le 4k+2$. Finally, for $1\le i\le m-1$,
\[(L_{i+1})_{4k+3, 4k+3}-(L_i)_{4k+3,4k+3}=+2 \mbox{ if } i\ne \lceil m/2\rceil.\]
Thus, $\{w(v_{i,4k+3})\;|\; 1\le i\le m\}$ consists of two different values.

\ms\nt Since each $L_i$ satisfies the condition of Lemma~\ref{lem-rowsum}, $\rs_j(L_i)$ is a strictly increasing function of $j$. Hence $\chi_{lat}(mK_n)\le n+1$.
\end{proof}

\begin{example} Take $n=7$ and $m=2$. So
\[\S_6=\begin{pmatrix}
+1 & -1 & +1 & -1 & +1 & -1\\
-1 & +1 & -1 & +1 & -1 & +1\\
+1 & -1 & -1 & +1 & +1 & -1\\
-1 & +1 & +1 & -1 & -1 & +1\\
+1 & -1 & +1 & -1 & +1 & -1\\
-1 & +1 & -1 & +1 & -1 & +1
\end{pmatrix} \rightarrow
\S_7=\left(\begin{array}{*{6}{c}|c}
+2 & -1 & +1 & -1 & +1 & -1 & -1\\
-1 & +2 & -1 & +1 & -1 & +1 & -1\\
+1 & -1 & -2 & +1 & +1 & -1 & +1\\
-1 & +1 & +1 & -2 & -1 & +1 & +1\\
+1 & -1 & +1 & -1 & +2 & -1 & -1\\
-1 & +1 & -1 & +1 & -1 & +2 & -1\\\hline
-1 & -1 & +1 & +1 & -1 & -1 & +2
\end{array}\right).
\]
\[\mathcal M=\left(\begin{array}{*{7}{c}}
* & -1 & +2 & -3 & +4 & -5 & -6\\
-1 & * & -7 & +8 & -9 & +10 & -11\\
+2 & -7 & * & +12 & +13 & -14 & +15\\
-3 & +8 & +12 & * & -16 & +17 & +18\\
+4 & -9 & +13 & -16 & * & -19 & -20\\
-5 & +10 & -14 & +17 & -19 & * & -21\\
-6 & -11 & +15 & +18 & -20 & -21 & *
\end{array}\right).
\]
\[L_1=\left(\begin{array}{*{7}{c}||c}
43 & 2 & 3 & 6 & 7 & 10 & 12 & 83\\
2 & 44 & 14 & 15 & 18 & 19 & 22 & 134\\
3 & 14 & 49 & 23 & 25 & 28 & 29 & 171\\
6 & 15 & 23 & 50 & 32 & 33 & 35 & 194\\
7 & 18 & 25 & 32 & 51 & 38 & 40 & 211\\
10 & 19 & 28 & 33 & 38 & 52 & 42 & 222\\
12 & 22 & 29 & 35 & 40 & 42 & 55 & 235
\end{array}\right),\quad
L_2=\left(\begin{array}{*{7}{c}||c}
45 & 1 & 4 & 5 & 8 & 9 & 11 & 83\\
1 & 46 & 13 & 16 & 17 & 20 & 21 & 134\\
4 & 13 & 47 & 24 & 26 & 27 & 30 & 171\\
5 & 16 & 24 & 48 & 31 & 34 & 36 & 194\\
8 & 17 & 26 & 31 & 53 & 37 & 39 & 211\\
9 & 20 & 27 & 34 & 37 & 54 & 41 & 222\\
11 & 21 & 30 & 36 & 39 & 41 & 56 & 234
\end{array}\right)
\]
So $\chi_{lat}(2K_7)\le 8$. \rsq
\end{example}

\begin{theorem}\label{thm-mKn-1}
For $m\ge 2$ and $n=4k+1\ge 5$, $\chi_{lat}(mK_n)\le \min\{n+3, n-1+m\}$.
\end{theorem}
\begin{proof} Similar to the proof of Theorem~\ref{thm-mKn-3} we will define a sign matrix $\mathcal S$ and a guide matrix $\mathcal M$ of order $4k+1$. We define a $(4k)\times (4k)$ matrix $\S_{4k}'$ first.

Let $\S'_4=\begin{pmatrix}
+1 & -1 & +1 & -1\\
-1 & +1 & -1 & +1\\
+1 & -1 & +1 & -1\\
-1 & +1 & -1 & +1\end{pmatrix}$ and
$\S'_{4k}=\left(\begin{array}{c|l}
\S'_{4k-4} & \begin{matrix}\S'_4 \\ \vdots \\ \S'_4\end{matrix}\\\hline
\begin{matrix}\S'_4 & \cdots & \S'_4\end{matrix} & \S_4
\end{array}\right)$ when $k\ge 2$.

\ms\nt Now, we define a symmetric sign matrix $\S_{4k+1}$ of order $4k+1$ using the same method as in the proof of Theorem~\ref{thm-mKn-3}. Note that the $(4k+1, 4k+1)$-entry of $\S_{4k+1}$ is $+4$.
\nt The definition of a guide matrix $\mathcal M$ and Stage~1 are similar to the proof of Theorem~\ref{thm-mKn-3}.

\ms\nt {\bf Stage 2: } Using a similar procedure to Stage~2 in the proof of Theorem~\ref{thm-mKn-3}, fill all diagonals of each $L_i$ except $(L_i)_{4k+1, 4k+1}$ by using $T(j)$ defined in the proof of Theorem~\ref{thm-mKn-3}, $1\le j\le 4k$.

\ms\nt Now we have to fill $[mN_4+4mk+1, mN_4+4mk+m]$ in the $(L_i)_{4k+1, 4k+1}$, here $N_4=\frac{n(n-1)}{2}$. Suppose $m=4s+m_0$, for some $s\ge 0$ and $0\le m_0<4$.

\ms\nt Suppose $s\ge 1$. Let $U(a)=\{mN_4+4km +a +4l\;|\; 0\le l\le s\}$ for $1\le a\le m_0$; and $U(a)=\{mN_4+4km +a +4l\;|\; 0\le l\le s-1\}$ for $m_0<a\le 4$. Let $U$ be the compound sequence $U^+(1)U^+(2)U^+(3)U^+(4)$. Let the $(4k+1,4k+1)$-entry of $L_i$ be the $i$-th term of $U$.

\nt Using a proof similar to Theorem~\ref{thm-mKn-3}, we have $\chi_{lat}(mK_n)\le n+3$.

\ms\nt Suppose $s=0$. That means $2\le m\le 4$. Then fill $[mN_4+4km+1, mN_4+4km+m]$ to the $(4k+1,4k+1)$-entry of $L_i$, respectively.

\ms\nt Using a proof similar to that above, $\rs_j(L_i)$ is a strictly increasing function of $j$ and hence we have $\chi_{lat}(mK_n)\le n-1+m$.

\ms\nt Combining these two cases, we conclude that $\chi_{lat}(mK_n)\le \min\{n+3, n-1+m\}$.
\end{proof}

\begin{corollary}\label{cor-A(mKnK1)-odd}
Let $m\ge 2$ and odd $n\ge 3$, $\chi_{lat}(A(mK_n, K_1))=n$.
\end{corollary}
\begin{proof} From the proofs of Theorems~\ref{thm-mKn-3} and \ref{thm-mKn-1}, we see that the $(n,n)$-th entries of all $L_i$'s are the largest $m$ labels. If we change $(L_i)_{n,n}$ to $N$ and denote this new matrix by $M_i$, then the matrix $M$ becomes a total labeling matrix of $A(mK_n, K_1)$ with $n$ difference row sums. Note that, $(L_i)_{n,n}$'s are identified to $A_{1,1}$. Thus $\chi_{lat}(A(mK_n, K_1))=n$.
\end{proof}

\begin{example} Let $n=9$, $m=2$. Then
\[\fontsize{7}{10}\selectfont \S_8'=\left(\begin{array}{*{4}{c}|*{4}{c}}
+1 & -1 & +1 & -1 & +1 & -1 & +1 & -1\\
-1 & +1 & -1 & +1 & -1 & +1 & -1 & +1\\
+1 & -1 & +1 & -1 & +1 & -1 & +1 & -1\\
-1 & +1 & -1 & +1 & -1 & +1 & -1 & +1\\\hline
+1 & -1 & +1 & -1 & +1 & -1 & +1 & -1\\
-1 & +1 & -1 & +1 & -1 & +1 & -1 & +1\\
+1 & -1 & +1 & -1 & +1 & -1 & -1 & +1\\
-1 & +1 & -1 & +1 & -1 & +1 & +1 & -1
\end{array}\right)\rightarrow
\S_9=\left(\begin{array}{*{4}{c}|*{4}{c}|c}
+2 & -1 & +1 & -1 & +1 & -1 & +1 & -1 & -1\\
-1 & +2 & -1 & +1 & -1 & +1 & -1 & +1 & -1\\
+1 & -1 & +2 & -1 & +1 & -1 & +1 & -1 & -1\\
-1 & +1 & -1 & +2 & -1 & +1 & -1 & +1 & -1\\\hline
+1 & -1 & +1 & -1 & +2 & -1 & +1 & -1 & -1\\
-1 & +1 & -1 & +1 & -1 & +2 & -1 & +1 & -1 \\
+1 & -1 & +1 & -1 & +1 & -1 & -2 & +1 & +1\\
-1 & +1 & -1 & +1 & -1 & +1 & +1 & -2 & +1\\\hline
-1 & -1 & -1 & -1 & -1 & -1 & +1 & +1 & +4
\end{array}\right)\]
Hence \[\M=\left(\begin{array}{*{4}{c}|*{4}{c}|c}
* & -1 & +2 & -3 & +4 & -5 & +6 & -7 & -8\\
-1 & * & -9 & +10 & -11 & +12 & -13 & +14 & -15\\
+2 & -9 & * & -16 & +17 & -18 & +19 & -20 & -21\\
-3 & +10 & -16 & * & -22 & +23 & -24 & +25 & -26\\\hline
+4 & -11 & +17 & -22 & * & -27 & +28 & -29 & -30\\
-5 & +12 & -18 & +23 & -27 & * & -31 & +32 & -33 \\
+6 & -13 & +19 & -24 & +28 & -31 & * & +34 & +35\\
-7 & +14 & -20 & +25 & -29 & +32 & +34 & * & +36\\\hline
-8 & -15 & -21 & -26 & -30 & -33 & +35 & +36 & *
\end{array}\right).\]
We have
\[L_1=\left(\begin{array}{*{9}{c}||c}
73 & 2 & 3	& 6	& 7	& 10 & 11 & 14 & 16 & 142\\
2 & 74 & 18 & 19 & 22 & 23 & 26	& 27 & 30 & 241\\
3 & 18 & 77 & 32 & 33 & 36 & 37	& 40 & 42 & 318\\
6 & 19 & 32 & 78 & 44 & 45 & 48	& 49 & 52 & 373\\
7 & 22 & 33	& 44 & 81 & 54 & 55	& 58 & 60 & 414\\
10 & 23	& 36 & 45 & 54 & 82& 62 & 63 & 66 & 441\\
11 & 26	& 37 & 48 & 55 & 62	& 87 & 67 & 69 & 462\\
14 & 27	& 40 & 49 & 58 & 63	& 67 & 88 & 71 & 477\\
16 & 30	& 42 & 52 & 60 & 66	& 69 & 71 & 89 & 495
\end{array}\right),
\]
\[L_2=\left(\begin{array}{*{9}{c}||c}
75 & 1 & 4 & 5 & 8 &9 & 12 & 13 & 15 & 142\\
1 &	76 & 17	& 20 & 21 & 24 & 25 & 28 & 29 & 241\\
4 &	17 & 79	& 31 & 34 & 35 & 38	& 39 & 41 & 318\\
5 &	20 & 31	& 80 & 43 &	46 & 47	& 50 & 51 & 373\\
8 & 21 & 34 & 43 & 83 & 53 & 56	& 57 & 59 & 414\\
9 & 24 & 35 & 46 & 53 &	84 & 61 & 64 & 65 & 441\\
12 & 25	& 38 & 47 &	56 & 61 & 85 & 68 &	70 & 462\\
13 & 28	& 39 & 50 & 57 & 64	& 68 & 86 & 72 & 477\\
15 & 29	& 41 & 51 & 59 & 65	& 70 & 72 & 90 & 492
\end{array}\right).
\]
Hence $\chi_{lat}(2K_9)\le 10$.  Clearly $\d(M)$ is larger than all vertex weights, so $\chi_{la}((2K_9)\vee K_1)\le 11$.

\ms\nt If we change $(L_2)_{9,9}$ to $89$ and identify the vertices $v_{1,9}$ with $v_{2,9}$, then we have a local antimagic total labeling for $A(2K_9, K_1)$ and hence $\chi_{lat}(A(2K_9, K_1))=9$. Clearly $\d(M)$ is larger than all vertex weights, so $\chi_{la}(A(2K_9, K_1)\vee K_1)=10$. \rsq
\end{example}

\begin{corollary}
Let $m\ge 2$ and odd $n\ge 3$, \[\chi_{la}(mK_{n+1})=\chi_{la}((mK_{n})\vee K_1)\le \begin{cases}
\min\{n+4, n+m\} & \mbox{ if }n\equiv 1\pmod{3},\\
n+2 & \mbox{ if }n\equiv 3\pmod{3}.\end{cases}\]
\end{corollary}
\begin{proof} From the proofs of Theorems~\ref{thm-mKn-3} and \ref{thm-mKn-1}, we see that the diagonals of $M$ are the largest $mn$ labels. Thus, $\mathcal D(M)>w(v_{m,n})$. So we have the corollary.
\end{proof}

\nt By Corollary~\ref{cor-A(mKnK1)-odd} and the same argument above, we have
\begin{corollary}
Let $m\ge 2$ and odd $n\ge 3$, $\chi_{la}(A(mK_{n+1},K_2))=\chi_{la}(A(mK_{n},K_1)\vee K_1)=n+1$.
\end{corollary}

\nt We have an ad hoc result for $A(mK_n, K_2)$ for $n$ odd as follows.

\begin{theorem}\label{thm-A(mKnK2)-3}
For $2\le m\le 3$ and odd $n\ge 3$, $\chi_{lat}(A(mK_n,K_2))=n$ and $\chi_{la}(A(mK_n, K_2)\vee K_1)=n+1$.
\end{theorem}

\begin{proof} The guide matrix is obtained from the guide matrix $\mathcal M$ of order $n$ defined in the proof of Theorem~\ref{thm-mKn-3} or Theorem~\ref{thm-mKn-1} by deleting the last two rows. Apply the same procedure as Stage~1 in the proof of Theorem~\ref{thm-mKn-3} or Theorem~\ref{thm-mKn-1}. It is equivalent to defining the matrix $\begin{pmatrix}L_i & B_i\end{pmatrix}$ of order $(n-2)\times n$ by deleting the last two rows of $L_i$ defined in the proof of Theorem~\ref{thm-mKn-3} or Theorem~\ref{thm-mKn-1}. Note that $(\mathcal M)_{n-1,n}=-\frac{n(n-1)}{2}$ is not used in this stage. So the labels used are $[1, mN_5]$, where $N_5=\frac{(n-2)(n+1)}{2}$.

\ms\nt{\bf Stage 2: } Similar to Stage~2 of the proof of Theorem~\ref{thm-mKn-3} or Theorem~\ref{thm-mKn-1}, we define the set $T(j)$ by labels $[mN_5+1, mN_5+m(n-2)+m]$ for $1\le j\le n-1$ and fill the diagonals of $L_i$'s. Note that, $mN_5+m(n-2)+m\le N$ and $T(n-1)$ is not used in this stage.

\nt By Lemma~\ref{lem-rowsum}, each $j$-th row sum of each $\begin{pmatrix}L_i & B_i\end{pmatrix}$ is a constant, and each $j$-th row sum of a fixed matrix $\begin{pmatrix}L_i & B_i\end{pmatrix}$ is a strictly increasing function of $j$, $1\le j\le n-2$.

\ms\nt{\bf Stage 3: } Use the remaining 3 labels to fill in the matrix $A_{1,2}, A_{1,1}, A_{2,2}$ in natural order. It is easy to see that the last two row sums of $M$ are distinct and larger than the other row sums of $M$.

\ms\nt It is easy to see that the diagonal entries of $M$ are the largest entries in the corresponding columns, and $\d(M)$ is greater than all weights of vertices.

\ms\nt Thus we have $\chi_{lat}(A(mK_n, K_2))=n$ and $\chi_{la}(A(mK_n, K_2)\vee K_1)=n+1$.
\end{proof}

\begin{example} Take $n=7$ and $m=2$.

\fontsize{8}{10}\selectfont
\[\begin{pmatrix}L_1 & B_1\end{pmatrix}=\left(\begin{array}{*{5}{c}|cc||c}
41 & 2 & 3 & 6 & 7 & 10 & 12 & 82\\
2 & 42 & 14 & 15 & 18 & 19 & 22 & 132\\
3 & 14 & 47 & 23 & 25 & 28 & 29 & 169\\
6 & 15 & 23 & 48 & 32 & 33 & 35 & 192\\
7 & 18 & 25 & 32 & 49 & 38 & 40 & 209\\\hline
10 & 19 & 28 & 33 & 38 & 52 & 50 & 230\\
12 & 22 & 29 & 35 & 40 & 50 & 53 & 241
\end{array}\right),\quad
\begin{pmatrix}L_2 & B_2\end{pmatrix}=\left(\begin{array}{*{5}{c}|cc||c}
43 & 1 & 4 & 5 & 8 & 9 & 11 & 81\\
1 & 44 & 13 & 16 & 17 & 20 & 21 & 132\\
4 & 13 & 45 & 24 & 26 & 27 & 30 & 169\\
5 & 16 & 24 & 46 & 31 & 34 & 36 & 192\\
8 & 17 & 26 & 31 & 51 & 37 & 39 & 209\\\hline
9 & 20 & 27 & 34 & 37 & 52 & 50 & 229\\
11 & 21 & 30 & 36 & 39 & 50 & 53 & 240
\end{array}\right).
\]
$w(u_6)=230+229-(50+52)=357$ and $w(u_7)=241+240-(50+53)=378$.
So $\chi_{lat}(A(2K_7,K_2))=7$. \rsq
\end{example}

\section{Conclusion and Open Problems}

\nt In this paper, $\chi_{lat}(A(mK_n,K_1))$ and $\chi_{la}(A(mK_n,K_1)\vee K_1)$ are determined except the following remaining case. 

\begin{problem} For $m\ge 2$ odd, $n\ge 3$, $n > r$, $r=0,2$, determine $\chi_{lat}(A(mK_n,K_r))$. \end{problem}

\nt From the motivation, we end this paper with the following problem.

\begin{problem} Determine $\chi_{lat}(A(mC_n,K_2))$ for $m\ge 2, n\ge 3$. \end{problem}







\begin{thebibliography}{99}
\bibitem{Arumugam} S. Arumugam, K. Premalatha, M. Bac\v{a} and A. Semani\v{c}ov\'{a}-Fe\v{n}ov\v{c}\'{i}kov\'{a}, Local antimagic vertex coloring of a graph, {\it Graphs Combin.}, {\bf33} (2017), 275--285.

\bibitem{Bensmail} J. Bensmail, M. Senhaji and K.S. Lyngsie, On a combination of the 1-2-3 conjecture and the antimagic labelling conjecture, {\it Discrete Math. Theoret. Comput. Sci.}, {\bf19(1)} (2017) \#22.

\bibitem{Bondy} J.A. Bondy, U.S.R. Murty, {\it Graph theory with applications}, New York, MacMillan, 1976.

%

\bibitem{Hagedorn} T.R. Hagedorn, Magic rectangles revisited, {\it Discrete Math.}, {\bf 207} (1999), 65--72.

\bibitem{Haslegrave} J. Haslegrave, Proof of a local antimagic conjecture, {\it Discrete Math. Theor. Comput. Sci.}, {\bf 20(1)} (2018), \#18.



\bibitem{Lau+LNS} G.C. Lau, J. Li, H.K. Ng and W.C. Shiu, Approaches which output infinitely many graphs with small local antimagic chromatic number, {\it Disc. Math. Algorithms Appl.}, (2022), accepted, arXiv:2009.01996.

\bibitem{LNS} G.C. Lau, H.K. Ng, and W.C. Shiu, Affirmative solutions on local antimagic chromatic number, {\it Graphs Combin.}, {\bf36} (2020), 1337--1354.


\bibitem{Lau+K+S} G.C. Lau, K. Schaeffer, W.C. Shiu, Every graph is local antimagic total and its applications, (2022), arxiv.org/abs/1906.10332.

\bibitem{Lau+SN} G.C. Lau, W.C. Shiu and H.K. Ng, On local antimagic chromatic number of graphs with cut-vertices, {\it Iran. J. Math. Sci. Inform.}, (2022) accepted, arXiv:1805.04801.

\bibitem{LSN1} G.C. Lau, W.C. Shiu and H.K. Ng, On local antimagic chromatic number of cycle-related join graphs, {\it Discuss. Math. Graph Theory}, {\bf 41} (2021) 133--152 DOI : 10.7151/dmgt.2177.

\bibitem{Lau+SS} G.C. Lau, W.C. Shiu and C.X. Soo, On local antimagic chromatic number of spider graphs, (2022) {\it J. Discrete Math. Sci. Cryptogr.},  DOI : 10.1080/09720529.2021.1892270

\bibitem{Premalatha+ALW} K. Premalatha, S. Arumugam, Y-C. Lee and T-M. Wang, Local antimagic chromatic number of trees - I, {\it J. Discrete Math. Sci. Cryptogr.}, (2020) DOI : 10.1080/09720529.2020.1772985.


\bibitem{Shiu+Lam+Lee2002} W.C. Shiu, P.C.B. Lam and S-M. Lee, On a construction of supermagic graphs, {\it JCMCC}, {\bf 42}, 147-160 (2002).
\bibitem{Zukerman} D. Zuckerman, Linear degree extractors and the inapproximability of max clique and chromatic number, {\it Theory Comput.}, {\bf3} (2007) 103--128, doi:10.4086/toc.2007.v003a006
\end{thebibliography}
\end{document}